\documentclass[a4paper,12pt,draft]{article}
\makeatletter
\@setfontsize\normalsize\@xiipt{15.5}%
\renewcommand*{\author}[1]{\gdef\@author{#1}\gdef\@pauthor{{\def\and{ --- }#1}}}
\renewcommand*{\title}[1]{\gdef\@title{#1}\gdef\@ptitle{#1}}
\if@twoside         
  \def\ps@draft{%
    \def\@oddfoot{\small\null\hfil\thepage\hfil}
    \let\@evenfoot\@oddfoot
    \def\@evenhead{\small\@date\hfil\slshape\@pauthor\hfil}
    \def\@oddhead{\small\null\hfil\slshape\@ptitle\hfil}
    \let\@mkboth\@gobbletwo
    \let\sectionmark\@gobble
    \let\subsectionmark\@gobble
   }
\else               
  \def\ps@draft{%
    \def\@oddfoot{\small\@date\hfil\slshape\@pauthor\hfil\upshape\thepage}
    \def\@oddhead{\small\null\hfil\slshape\@ptitle\hfil}
    \let\@mkboth\@gobbletwo
    \let\sectionmark\@gobble
    \let\subsectionmark\@gobble
  }
\fi
\newcommand*{\keywords}[1]{\gdef\@keywords{#1}}
\keywords{}
\newcommand{\keywordsname}{Key words and phrases}
\newcommand*{\subjclass}[1]{\gdef\@subjclass{#1}}
\subjclass{}
\newcommand{\subjclassname}{1991 AMS Mathematics Subject Classification}
\def\@maketitle{%
  \newpage
  \null
  \vskip 2em%
  \begin{center}%
    {\Large\bfseries \@title \par}%
    \vskip 1.5em%
    {\small\scshape
      \lineskip .5em%
      \begin{tabular}[t]{c}%
        \@author
      \end{tabular}\par
    }%
    \vskip 1em%
    {\small\@date}
  \end{center}%
  \par
  \vskip 1.5em
  \begingroup
    \let\@makefnmark\relax \let\@thefnmark\relax
    \ifx\@empty\@subjclass\else
       \@footnotetext{{\itshape\subjclassname}.\enspace\@subjclass.}
    \fi
    \ifx\@empty\@keywords\else
       \@footnotetext{{\itshape\keywordsname}.\enspace\@keywords.}
    \fi
  \endgroup
}
\renewcommand{\part}{\par
   \addvspace{4ex}%
   \@afterindentfalse
   \secdef\@part\@spart}
\def\@part[#1]#2{%
    \ifnum \c@secnumdepth >\m@ne
      \refstepcounter{part}%
      \addcontentsline{toc}{part}{\thepart\hspace{1em}#1}%
    \else
      \addcontentsline{toc}{part}{#1}%
    \fi
    {\parindent \z@ \raggedright
     \interlinepenalty \@M
     \reset@font
     \ifnum \c@secnumdepth >\m@ne
       \large\bfseries \partname~\thepart
       \par\nobreak
     \fi
     \Large \bfseries #2%
     \markboth{}{}\par}%
    \nobreak
    \vskip 3ex
    \@afterheading}
\def\@spart#1{%
    {\parindent \z@ \raggedright
     \interlinepenalty \@M
     \reset@font
     \Large \bfseries #1\par}%
     \nobreak
     \vskip 3ex
     \@afterheading}
\def\@endpart{\vfil\newpage
              \if@twoside
                \hbox{}%
                \thispagestyle{empty}%
                \newpage
              \fi
              \if@tempswa
                \twocolumn
              \fi}
\renewcommand{\section}{\@startsection {section}{1}{\z@}%
                                   {-3.5ex \@plus -1ex \@minus -.2ex}%
                                   {2.3ex \@plus.2ex}%
                                   {\reset@font\large\bfseries}}
\renewcommand{\subsection}{\@startsection{subsection}{2}{\z@}%
                                     {-3.25ex\@plus -1ex \@minus -.2ex}%
                                     {1.5ex \@plus .2ex}%
                                     {\reset@font\normalsize\bfseries}}
\renewcommand{\subsubsection}{\@startsection{subsubsection}{3}{\z@}%
                                     {-3.25ex\@plus -1ex \@minus -.2ex}%
                                     {1.5ex \@plus .2ex}%
                                     {\reset@font\normalsize\bfseries}}
\renewcommand{\paragraph}{\@startsection{paragraph}{4}{\z@}%
                                    {3.25ex \@plus1ex \@minus.2ex}%
                                    {-1em}%
                                    {\reset@font\normalsize\bfseries}}
\renewcommand{\subparagraph}{\@startsection{subparagraph}{5}{\parindent}%
                                       {3.25ex \@plus1ex \@minus .2ex}%
                                       {-1em}%
                                      {\reset@font\normalsize\bfseries}}
\renewcommand{\theenumi}{\alph{enumi}}

\renewcommand{\theenumii}{\roman{enumii}}

\renewcommand{\p@enumii}{\theenumi.}
\renewcommand{\theenumiii}{\Alph{enumiii}}

\renewcommand{\p@enumiii}{\theenumi.\theenumii.}

\renewcommand{\p@enumiv}{\p@enumiii\theenumiii.}
\@addtoreset{equation}{section}

\RequirePackage{amsthm}
\renewenvironment{proof}[1][\proofname]{\par
  \normalfont
  \topsep6\p@\@plus6\p@ \trivlist
  \item[\hskip\labelsep\slshape
    #1\@addpunct{.}]\ignorespaces
}{%
  \qed\endtrivlist
}
\theoremstyle{plain}
\newtheorem{theorem}{Theorem}[section]
\newtheorem{proposition}[theorem]{Proposition}
\newtheorem{lemma}[theorem]{Lemma}
\newtheorem{corollary}[theorem]{Corollary}

\theoremstyle{definition}

\newtheorem{example}[theorem]{Example}

\newtheorem{remark}[theorem]{Remark}

\RequirePackage{amsmath}
\RequirePackage{amssymb}
\renewcommand{\emptyset}{\varnothing}

\newcommand{\C}{\mathbb{C}}
\newcommand{\N}{\mathbb{N}}
\newcommand{\R}{\mathbb{R}}

\DeclareMathOperator{\id}{id}

\newcommand{\plim}[1][]{\mathop{\varprojlim}\limits_{#1}}
\DeclareMathOperator{\coker}{coker}

\renewcommand{\to}[1][]{\xrightarrow[#1]{}}

\newcommand{\Endo}[1][]{\mathrm{End}_{\raise1.5ex\hbox to.1em{}#1}}
\newcommand{\Hom}[1][]{\mathrm{Hom}_{\raise1.5ex\hbox to.1em{}#1}}
\newcommand{\RHom}[1][]{\mathrm{RHom}_{\raise1.5ex\hbox to.1em{}#1}}
\newcommand{\Ext}[2][]{\mathrm{Ext}_{\raise1.5ex\hbox to.1em{}#1}^{#2}}
\newcommand{\THom}[1][]{\mathrm{THom}_{\raise1.5ex\hbox to.1em{}#1}}

\newcommand{\Tens}[1][]{\mathbin{\otimes_{\raise1.5ex\hbox to-.1em{}#1}}}
\newcommand{\LTens}[1][]{\mathbin{\otimes_{\raise1.5ex\hbox to-.1em{}#1}^{L}}}
\newcommand{\Tor}[2][]{\mathrm{Tor}^{\raise1.5ex\hbox to.1em{}#1}_{#2}}

\def\sha{\mathcal{A}}
\def\shb{\mathcal{B}}
\def\shc{\mathcal{C}}
\def\shd{\mathcal{D}}
\def\she{\mathcal{E}}

\def\shi{\mathcal{I}}

\def\sho{\mathcal{O}}

\def\sht{\mathcal{T}}
\def\shu{\mathcal{U}}

\newcommand{\sect}{\Gamma}

\renewcommand{\hom}[1][]{{\mathcal{H}om}_{\raise1.5ex\hbox to.1em{}#1}}
\newcommand{\rhom}[1][]{{R\mathcal{H}om}_{\raise1.5ex\hbox to.1em{}#1}}
\newcommand{\ext}[2][]{{\mathcal{E}xt}_{\raise1.5ex\hbox to.1em{}#1}^{#2}}
\newcommand{\thom}[1][]{{T\mathcal{H}om}_{\raise1.5ex\hbox to.1em{}#1}}

\newcommand{\tens}[1][]{\mathbin{\otimes_{\raise1.5ex\hbox to-.1em{}#1}}}
\newcommand{\ltens}[1][]{\mathbin{\otimes_{\raise1.5ex\hbox to-.1em{}#1}^{L}}}
\newcommand{\tor}[2][]{{\mathcal{T}or}^{\raise1.5ex\hbox to.1em{}#1}_{#2}}
\newcommand{\wtens}{\mathbin{\mathop{\otimes}\limits^{{}_{\mathrm{w}}}}}
\newcommand\etens{\mathbin{\boxtimes}}
\DeclareMathOperator{\supp}{supp}
\newcommand{\oim}[1]{{#1}_*}
\newcommand{\eim}[1]{{#1}_!}

\newcommand{\opb}[1]{#1^{-1}}

\newcommand{\GHom}[1][]{\mathrm{GHom}_{\raise1.5ex\hbox to.1em{}#1}}
\newcommand{\GExt}[2][]{\mathrm{GExt}_{\raise1.5ex\hbox to.1em{}#1}^{#2}}
\newcommand{\FHom}[1][]{\mathrm{FHom}_{\raise1.5ex\hbox to.1em{}#1}}
\newcommand{\ghom}[1][]{{\mathcal{GH}om}_{\raise1.5ex\hbox to.1em{}#1}}
\newcommand{\gext}[2][]{{\mathcal{GE}xt}_{\raise1.5ex\hbox to.1em{}#1}^{#2}}
\newcommand{\fhom}[1][]{{\mathcal{FH}om}_{\raise1.5ex\hbox to.1em{}#1}}

\newcommand{\tenstop}[1][]{\mathbin{\hat{\otimes}_{\raise1.5ex\hbox to-.1em{}#1}}}
\newcommand{\homtop}[1][]{\mathcal{L}_{\raise1.5ex\hbox to.1em{}#1}}
\newcommand{\Homtop}[1][]{\mathrm{L}_{\raise1.5ex\hbox to.1em{}#1}}

\newcommand{\A}{\mathcal{A}}

\def\absdoim#1{\underline{#1}_*}
\def\reldoim[#1]#2{\underline{#2}_{|{#1}*}}
\def\doim{\@ifnextchar [{\reldoim}{\absdoim}}
\def\absdeim#1{\underline{#1}_*}
\def\reldeim[#1]#2{\underline{#2}_{|{#1}*}}
\def\deim{\@ifnextchar [{\reldeim}{\absdeim}}
\def\absdopb#1{\underline{#1}^{-1}}
\def\reldopb[#1]#2{\underline{#2}_{|{#1}}^{-1}}
\def\dopb{\@ifnextchar [{\reldopb}{\absdopb}}
\def\absboim#1{\underline{\underline{#1}}_*}
\def\relboim[#1]#2{\underline{\underline{#2}}_{|{#1}*}}
\def\boim{\@ifnextchar [{\relboim}{\absboim}}
\def\absbeim#1{\underline{\underline{#1}}_*}
\def\relbeim[#1]#2{\underline{\underline{#2}}_{|{#1}*}}
\def\beim{\@ifnextchar [{\relbeim}{\absbeim}}
\def\absbopb#1{\underline{\underline{#1}}^*}
\def\relbopb[#1]#2{\underline{\underline{#2}}_{|{#1}}^*}
\def\bopb{\@ifnextchar [{\relbopb}{\absbopb}}

\RequirePackage{xypic}

\def\etens{\mathbin{\boxtimes}}
\newcount\temp
\temp=\catcode`\@
\catcode`\@=11
\def\@bletens{\mathbin{\etens^{L}}}
\def\@letens_#1{\mathbin{\etens_{\raise1.5ex\hbox to-.1em{}#1}^{L}}}
\def\letens{\@ifnextchar _{\@letens}{\@bletens}}
\catcode`\@=\temp

\newcommand{\ihom}[1][]{{\shi hom}_{\raise1.5ex\hbox to.1em{}#1}}
\newcommand{\eeim}[1]{{#1}_{!!}}

\def\phi{{\varphi}}

\newcommand{\indlim}[1][]{\mathop{\varinjlim}\limits_{#1}}
\newcommand{\prolim}[1][]{\mathop{\varprojlim}\limits_{#1}}

\def\epsilon{\varepsilon}

\newcommand{\indc}{{\rm Ind}(\shc)}
\newcommand{\indcc}{{\rm Ind}(\shc')}

\newcommand{\inddlim}[1][]{{``{\indlim[#1]}"}}

\def\epito{{\twoheadrightarrow}}
\newcommand{\ssubset}{\subset\subset}

\newcommand{\md[1]}{{\rm Mod}({#1})}

\newcommand{\mdc[1]}{{\on{Mod}}^c({#1})}

\newcommand{\II[1]}{{\rm I}({#1})}
\def\isha{{\II[\sha]}}

\def\ishb{{\II[\shb]}}

\newcommand{\eqn}{\begin{eqnarray*}}
\newcommand{\eneqn}{\end{eqnarray*}}
\newcommand{\eq}{\begin{eqnarray}}
\newcommand{\eneq}{\end{eqnarray}}
\newenvironment{nnum}{
  \begin{enumerate}
   \itemsep=0pt

   }
  {\end{enumerate}}
\newenvironment{anum}{
  \begin{enumerate}
     \itemsep=0pt

   }
  {\end{enumerate}}
\newcommand{\bnum}{\begin{nnum}}
\newcommand{\enum}{\end{nnum}}
\newcommand{\banum}{\begin{anum}}
\newcommand{\eanum}{\end{anum}}

\newcommand{\Lemma}{\begin{lemma}}
\newcommand{\enlemma}{\end{lemma}}
\newcommand{\Cor}{\begin{corollary}}
\newcommand{\encor}{\end{corollary}}
\newcommand{\Prop}{\begin{proposition}}
\newcommand{\enprop}{\end{proposition}}

\def\etens{\mathbin{\boxtimes}}
\newcount\temp
\temp=\catcode`\@
\catcode`\@=11
\def\@bletens{\mathbin{\etens^{L}}}
\def\@letens_#1{\mathbin{\etens_{\raise1.5ex\hbox to-.1em{}#1}^{L}}}
\def\letens{\@ifnextchar _{\@letens}{\@bletens}}
\catcode`\@=\temp

\def\rop{{\rm op}}

\def\phi{{\varphi}}

\newcommand{\rc[1]}{{\Brc}({#1})}
\newcommand{\rcc[1]}{{\Brc}^c({#1})}
\def\Irckx{{{\rm I}_{\R\hbox{-}c}(k_X)}}
\def\wtens{\stackrel{{\rm w}}{\tens}} 
\makeatother

\begin{document}
\author{Masaki Kashiwara \& Pierre Schapira}
\title{Ind-Sheaves, distributions, and microlocalization}
\date{04/1999\footnote{Talk given at Ecole Polytechnique, May 1999}}
\maketitle

\section{Introduction}
If $\shc$ is an abelian category, the category $\indc$ of ind-objects
of $\shc$ has many remarkable properties: 
it is much bigger than $\shc$, it
contains $\shc$, and furthermore it is dual (in a certain sense) to
$\shc$.

We introduce here the category of ind-sheaves on a locally compact 
space $X$ as the category of
ind-objects of the category of sheaves with compact supports. 
This construction has some analogy with that of distributions: the
space of distributions is bigger than that of functions, and is dual
to that of functions with compact support. This last condition implies
the local nature of distributions, and similarly, one proves that
the category of ind-sheaves defines a stack.

In our opinion, what makes the theory of ind-sheaves on manifolds really
interesting is twofold.

(i) Ind-sheaves allows us to treat in the formalism of sheaves (the
``six operations'') functions with growth conditions.
For example, on a complex manifold $X$, one can define the ind-sheaf of 
``tempered holomorphic functions'' $\sho_X^t$ or the ind-sheaf 
of ``Whitney holomorphic functions'' $\sho_X^{\rm w}$, and obtain for
example the sheaf of Schwartz's distributions 
using Sato's construction of hyperfunctions, 
simply replacing $\sho_X$ with $\sho_X^t$.

(ii) On a real manifold, one can construct a microlocalization functor $\mu_X$
which sends sheaves
(i.e. objects of the derived category of sheaves) on $X$ to
ind-sheaves on $T^*X$, and the Sato
functor of microlocalization along a submanifold
$M\subset X$ (see \cite{S-K-K}) becomes the usual functor 
$\hom(\opb{\pi}(k_M),\cdot)$ (where $\pi:T^*X\to X$ is the projection)
composed with $\mu_X(\cdot)$.

When combining (i) and (ii), one can treat in a unified way various
objects of classical analysis. 

The results presented here are extracted from \cite{K-S3}.
We refer to \cite{K-S1} for an exposition on derived categories and sheaves,
and to \cite{Gr-V} for the theory of ind-objects.

\section{Ind-objects}
Let $\shu$ be a universe (see \cite{McL}).
If $\shc$ is a $\shu$-category, 
one denotes by $\shc^\vee$ the category
of contravariant functors from $\shc$ to {\bf Set},
the category of ($\shu$-small) sets. One sends $\shc$ into $\shc^\vee$
by the fully faithful functor $h^\vee: X\mapsto\Hom[\shc](\cdot,X)$.

Let $I$ be a small filtrant category and let 
$i\mapsto X_i$ be an inductive system in $\shc$ indexed by $I$.
One denotes by $\inddlim[i]X_i$ the object of $\shc^\vee$ defined by
$$
\shc\ni Y\mapsto \indlim[i\in I]\Hom[\shc](Y,X_i).
$$
An ind-object in $\shc$ is an object of $\shc^\vee$
which is isomorphic to such a $\inddlim[i]X_i$.
One denotes by $\indc$ the full subcategory of 
$\shc^\vee$ consisting of ind-objects and identifies
$\shc$ with a full subcategory of $\indc$ by the functor $h^\vee$.

From now on we shall assume that $\shc$ is abelian. 
The subcategory $\shc^{\vee,add}$ of $\shc^\vee$
consisting of additive functors is clearly abelian and the functor
$h^\vee:\shc\to\shc^{\vee,add}$ is left exact. 

Since the functor $\Hom$ is left exact and filtrant inductive limits
are exact, any ind-object defines a left
exact functor on $\shc$. The converse holds if $\shc$ is small.

\begin{example} Let $k$ be a field. Define
the inductive system in $\md[k]$ indexed by $\N$ by setting $X_n=k^n$,
with the natural injections $X_n\hookrightarrow X_{n+1}$. Then 
$\inddlim[n]X_n$ is the functor
$Y\mapsto\indlim[n]\Hom[k](Y,X_n)$. 
This object does not belong to $\md[k]$ since there are no vector
space $Z$ such that $\indlim[n]\Hom[k](Y,X_n)\simeq \Hom[k](Y,Z)$ for
all vector spaces $Y$.
\end{example}

If $f:X\to Y$ is a morphism in $\indc$, one can construct a small
filtrant category $I$ such that $X=\inddlim[i]X_i$, $Y=\inddlim[i]Y_i$
and $f=\inddlim[i]f_i$, with $f_i:X_i\to Y_i$. Then $\inddlim[i]\ker f_i$
and $\inddlim[i]\coker f_i$ are the kernel and the cokernel of $f$ in $\indc$.

\begin{theorem}
\bnum
\item The category $\indc$ is abelian.
\item The natural functor $\shc\to \indc$ is fully faithful and
exact, and $\shc$ is thick in $\indc$.
\item The natural functor $\indc\to \shc^{\vee,add}$ is fully faithful
and left exact. 
\item The category $\indc$ admits right exact small inductive limits. 
Moreover, inductive limits over small filtrant categories are exact.
\item If small products of objects of $\shc$ exist in $\indc$, then
  $\indc$ admits small projective limits and such limits are left exact.
\enum
\end{theorem}
We shall denote by $\inddlim$ the usual inductive limit in
$\indc$. With this convention, one can identify $\shc$ with its
image in $\indc$ by $h^\vee$ without risk of confusion.

Consider an additive functor $F:\shc\to\shc'$. It defines an additive functor
$IF:\indc\to\indcc$. If there is no
risk of confusion, we shall write $F$ instead of $IF$. If $F$ is right exact,
or left exact, or fully faithhful, so is $IF$. 

We shall have to consider the derived category of $\indc$.
One proves that the natural functor $D^b(\shc)\to D_{\shc}^b (\indc)$
is an equivalence.

Even if $\shc$ has enough injectives, we cannot prove that
$\indc$ has the same property. However, quasi-injective objects,
i.e. ind-objects which are exact functors on $\shc$, are sufficient
for many applications.

\section{Ind-sheaves}
Let $X$ be a locally compact topological space which is countable at
infinity. Let $k$ be a field and let $\sha$ be a sheaf of commutative 
$k$-unitary algebras on $X$.

One denotes by
$\md[\sha]$ the abelian category of sheaves of $\sha$-modules, and by
$\mdc[\sha]$ the full subcategory consisting of sheaves with compact
support.

We call an object of ${\rm Ind}(\mdc[\A])$ an ind-sheaf on $X$ and we
set for short:
$$\II[\A]:= {\rm Ind}(\mdc[\A]).$$

If $U$ is an open subset of $X$, one defines the restriction
of $F$ to $U$ as follows. If $F=\inddlim[i] F_i$, one
sets
$$F\vert_U=\inddlim[i,V\ssubset U](F_{iV}\vert_U).$$
Then one proves that for $F$ and $G$ in $\II[\A]$, the presheaf 
$U\mapsto \Hom[{\II[\A\vert_U]}](F\vert_U,G\vert_U)$ is
a sheaf. We shall denote it by $\hom(F,G)$.

In fact, there is a better result: $U\mapsto\II[\A\vert_U]$ is a
stack. Roughly speaking, this means a sheaf of categories.

Note that ${\rm Ind}(\md[\A])$ is not a stack.
\begin{example}
Let $X=\R$, and let $\sha=k_X$.
Let $F=k_X, G_n=k_{[n,+\infty[}, G=\inddlim[n]G_n$.
Then $G\vert_U=0$ for any relatively compact open subset $U$ of $X$.
On the other hand,
$\Hom[{{\rm Ind}(\md[k_X])}](k_X,G)\simeq\indlim[n]\Hom[k_X](k_X,G_n)\simeq k$.
\end{example}

We construct the functors
\begin{eqnarray*}
\iota_X:\md[\A]\to \II[\A],&&
\iota_XF= \inddlim[U\subset\subset X]F_U,\\
\alpha_X:\II[\A]\to\md[\A],&&
\alpha_X(\inddlim[i]F_i)=\indlim[i]F_i.
\end{eqnarray*}
The functor $\alpha_X$ admits a left adjoint 
$\beta_X:\md[\A]\to\II[\A]$. 

These functors satisfy: 
\bnum
\item
$\iota_X$ and $\alpha_X$ are exact, fully faithful, 
and commute with $\indlim$ and $\prolim$, 
\item $\beta_X$ is right exact, fully faithful and commutes with $\indlim$, 
and if $\sha$ is left coherent, $\beta_X$ is exact.
\item $\alpha_X$ is left adjoint to $\iota_X$ and is right adjoint to 
$\beta_X$,
\item $\alpha_X\circ\iota_X\simeq \id$ and $\alpha_X\circ\beta_X\simeq \id$.
\enum
If $Z$ is locally closed in $X$ and $F$ is a sheaf on $X$, recall
that the sheaf $F_Z$ is $0$ on $X\setminus Z$ and is isomorphic to
$F\vert_Z$ on $Z$. 
We set
$$\widetilde{\sha_Z}=\beta_X(\sha_Z).$$
Since $\beta_X$ is right exact and commutes with $\indlim$, it is
characterized by its values on the sheaves $\sha_U$, $U$ open. 
If $U$ is open and $S$ is closed, one has
\begin{eqnarray*}
\widetilde{\sha_U}\simeq \inddlim[V\ssubset U]\sha_V,
\ \mbox{ $V$ open, and  }\quad&&
\widetilde{\sha_S}\simeq \inddlim[V\supset S]\sha_{\bar{V}},\,\ V\mbox{ open}.
\end{eqnarray*}
Note that $\widetilde{\sha_U}\to \sha_U$ is a monomorphism and 
$\widetilde{\sha_S}\to \sha_S$ is an epimorphism.

\begin{example}\label{ex:indshzeroonU}
Let $X$ be a real manifold of dimension $n\geq 1$ and let $\sha=k_X$.
For short, we shall write $k_Z$ instead of $(k_X)_Z$. Let $a\in X$.
Define $N_a\in\II[k_X]$ by the exact sequence
\begin{equation}\label{eq:Nx}
0\to N_a\to\widetilde{k_{\{a\}}}\to k_{\{a\}}\to 0.
\end{equation}
The derived functor $\RHom[{\II[k_X]}](\cdot,\cdot)$ is 
well-defined and moreover if $G$ and $F_i$ are sheaves on $X$, 
\begin{equation}\label{eq:rhomlim1}
H^k\RHom[{\II[k_X]}](G,\inddlim[i]F_i)\simeq\indlim[i]H^k\RHom[k_X](G,F_i).
\end{equation}
Since $\indlim[a\in V]H^n_{\{a\}}(X;k_{\overline{V}})\neq 0$, we find that the morphism $\widetilde{k_{\{a\}}}\to k_{\{a\}}$ is not an isomorphism,
hence $N_a\neq 0$.

On the other hand, for any open neighborhood $U$ of $a$,
we have
$$\RHom[{\II[k_X]}](k_U,N_a)\simeq 0.$$ 
\end{example}

\section{Operations on ind-sheaves} 

We define the functors internal tensor product, denoted $\tens$, and
internal $\hom$, denoted $\ihom$
\begin{eqnarray*} 
\tens: \II[\sha]\times\II[\sha] \to \II[\sha],\\
\ihom: \II[\sha]^\rop\times\II[\sha]\to \II[\sha],
\end{eqnarray*} 
by the formulas:
\begin{eqnarray*} 
\inddlim[i]F_i\tens\inddlim[j]G_j
                          &=&\inddlim[i,j](F_i\tens G_j),\\
\ihom(\inddlim[i]F_i,\inddlim[j]G_j)
                          &=&\prolim[i]\inddlim[j]\hom(F_i,G_j).
\end{eqnarray*} 
One has:
\begin{eqnarray*} 
\alpha_X\ihom(F,G)\simeq \hom(F,G).
\end{eqnarray*} 
The functor $\tens$ is right exact and commutes with $\inddlim$
and the functor $\ihom$ is left exact. 
The functors $\tens$ and $\ihom$ are adjoint:
\begin{eqnarray*}
\Hom[{\II[\sha]}](F\tens[\sha] K,G)
       \simeq \Hom[{\II[\sha]}](G,\ihom[\sha](K,G)),\\
\end{eqnarray*}

Now consider a continuous map $f:X\to Y$ of locally compact spaces,
and assume $\sha=\opb{f}\shb$, for a sheaf of rings $\shb$ on $Y$.

If $G=\inddlim[i]G_i$ is an ind-sheaf on $Y$, we define
$$
\opb{f}G=\inddlim[i](\opb{f}G_i)_U,\, U\ssubset X,\,\mbox{ $U$ open}.
$$
The functor $\opb{f}:\II[\shb]\to\isha$ is exact and commutes with
$\inddlim$ and $\tens\,$. Moreover, it commutes with the functors
$\iota,\alpha$ and $\beta$.

If $F=\inddlim[i]F_i$ is an ind-sheaf on $X$, we define:
$$
\oim{f}F=\prolim[K]\inddlim[i]\oim{f}F_{iK},\mbox{ $K$ compact}.
$$
The functor $\oim{f}:\isha\to\ishb$ is left  exact and commutes with
$\plim$. Moreover, it commutes with the functors $\iota$ and $\alpha$.

The two functors $\opb{f}$ and $\oim{f}$ are adjoint. More precisely,
for $F\in\isha$ and $G\in\ishb$, we have
$$\Hom[{\isha}](\opb{f}G,F)\simeq \Hom[{\ishb}](G,\oim{f}F).$$

If $F=\inddlim[i]F_i$ is an ind-sheaf on $X$, we define  $\eeim{f}F$ by
the formula
$$\eeim{f}\inddlim[i]F_i=\inddlim[i]\eim{f}F_i.$$
Note that the natural morphism $\eeim{f}\iota_XF\to
\iota_Y\eim{f}F$ is not an isomorphism in general.
The functor $\eeim{f}$ is left exact and commutes with $\inddlim$.
Moreover, it commutes with the functor $\alpha$.

There is a base change formula as well as a
projection formula for ind-sheaves.
\begin{example}
For $F\in \II[k_X]$ and $x\in X$ set $F_x=\opb{j_x}F$, where
$j_x:\{x\}\hookrightarrow X$. 
Let $N_a$ be as in Example \ref{ex:indshzeroonU}. Then $(N_a)_x\simeq
0$ for all $x\in X$. On the other hand, one shows that the functor 
$F\mapsto \prod_{x\in X}(F\tens \widetilde{k_{\{x\}}})$ is faithful.
\end{example}

\section{Construction of ind-sheaves}

We assume that $X$ is  a real analytic manifold
and denote by $\sht$ the family of open relatively compact
subanalytic subsets of $X$.
We denote by $\rc[k_X]$ the (small) category of $\R$-constructible sheaves of
$k$-vector spaces on $X$ and by $\rcc[k_X]$ the full subcategory
consisting of $\R$-constructible sheaves with compact support. 
Note that for any
$F\in\rcc[k_X]$ there is an exact sequence $F^1\to F^0\to F\to 0$,
with $F^0$ and $F^1$ finite direct sums of sheaves $k_U$ with $U\in\sht$.

We set
$$\Irckx={\rm Ind}(\rcc[k_X]).$$
The natural functor $\iota:\rcc[k_X]\to \mdc[k_X]$ defines the fully
faithful exact functor
\begin{equation}\label{functiota}
\iota:\Irckx\to\II[k_X].
\end{equation}
On the other hand, the forgetful functor 
$\mdc[k_X]^{\vee,add}\to \rcc[k_X]^{\vee,add}$ induces a functor
\begin{equation}\label{functrho}
\rho:\II[k_X]\to\Irckx.
\end{equation}
The next theorem
formulates a previous result of 
\cite{K-S2} in the language of ind-sheaves.

\begin{theorem}\label{th:constindsh1}
Let $F$ be a presheaf of $k$-vector spaces on $\sht$. Assume:
\bnum
\item $F(\emptyset)=0$,
\item for any $U$ and $V$ in $\sht$, the sequence 
$0\to F(U\cup V)\to F(U)\oplus F(V)\to F(U\cap V)$
is exact.
\enum
Then 
there exists $F^+_\sht\in\Irckx$ such that 
for any $U\in\sht$,
\begin{equation}
\Hom[{\Irckx}](k_U,F^+_\sht)\simeq F(U).
\end{equation}
\end{theorem}
We denote by $F^+$ the image of $F^+_\sht$ in $\II[k_X]$ by the 
functor $\iota$ of (\ref{functiota}).
If $\shd$ is a sheaf of (not necessarily commutative) $k$-algebras
and $F$ a presheaf of $\shd$-modules, then $F^+$ belongs to 
$\md[\shd,{\II[k_X]}]$, the subcategory of $\II[k_X]$ of ind-sheaves
endowed with a $\shd$-action. Recall that $F\in \md[\shd,{\II[k_X]}]$
means that $F$ is endowed with a morphism of sheaves of unitary rings
$\shd\to\she nd_{\II[k_X]}(F)$.

\section{Some classical ind-sheaves}\label{sect:cinftytemp}

In this section, the base field $k$ is $\C$. 
We denote by 
$\shc_X^\infty$ the sheaf of complex valued functions
of class $\shc^\infty$, 
and by $\shd_X$ the sheaf of real analytic finite-order differential operators.

(a) Let $f\in \shc^{\infty}(U)$. One says that $f$ has
{\it  polynomial growth} at $p\in X$ if for a local coordinate system
$(x_1,\cdots,x_n)$ around $p$, there exist a sufficiently small
compact neighborhood $K$ of $p$ and a positive integer $N$ such that
\begin{eqnarray}
&\sup_{x\in K\cap U}\big( dist(x,X\setminus U)\big)^N\vert f(x)\vert
<\infty\,.&
\end{eqnarray}
We say that $f$ is {\it tempered at $p$} if all of its derivatives are
of polynomial growth at $p$.
We say that $f$ is {\it tempered} if it is tempered at any point of $X$.
One denotes by $\shc^{\infty,t}(U)$ the $\C$-vector subspace of
$\shc^{\infty}(U)$ consisting of tempered functions.
By a theorem of Lojaciewicz, 
the contravariant functor $\shc^{\infty,t}(\cdot)$ defined on the
category $\sht$ of open relatively compact subanalytic sets is
exact. Hence, we may apply Theorem \ref{th:constindsh1}, 
and we get an ind-sheaf:
$$\shc^{\infty,t}_X\in\md[\shd_X,{\II[\C_X]}].$$

(b) If $S$ is 
closed in $X$, one denotes by $\shi_S ^\infty(X)$ the ideal of 
$\shc^{\infty}(X)$ consisting of functions which vanish on $S$
with infinite order, and for an  open subset $U$ in $X$, we set
$\sect(X;\C_U\wtens\shc^{\infty}_X)=\shi_{X\setminus U}^\infty(X)$.
Again by a theorem
of Lojaciewicz, the (covariant) functor 
$\sect(X;\cdot\wtens\shc^{\infty}_X)$
defined on the
category $\sht$ is exact. It extends to an exact functor 
$\sect(X;\cdot\wtens\shc_X^{\infty})$ on the
category $\rcc[\C_X]$ and we may define 
$$\shc^{\infty,{\rm w}}_X\in\md[\shd_X,{\II[{\rc[\C_X]}]}]$$
by the formula
$$\shc^{\infty,{\rm w}}_X(F)=\sect(X;H^0(D'F)\wtens\shc_X^{\infty})$$
where $F\in \rcc[\C_X]$ and $D'F=\rhom[\C_X](F,\C_X)$.

(c) Similarly, replacing $\wtens$ by $\tens$ in the above formula,
we find an ind-sheaf $\shc^{\infty,\omega}_X$ which, in fact, is
isomorphic to $\beta_X\shc^\infty_X$.

(d) Now assume that $X$ is a complex manifold with structure sheaf
$\sho_X$, denote by $\overline{X}$ the complex conjugate 
manifold and by $X^\R$ the
underlying real manifold, identified with the diagonal of
$X\times\overline{X}$. 
We define the objects $\sho_X^t$, $\sho_X^{\rm w}$ and by $\sho_X^\omega$
considering the Dolbeault complexes with coefficients in
$\shc^{\infty,t}_X$, $\shc^{\infty,{\rm w}}_X$ and $\shc^{\infty,\omega}_X$:
\begin{eqnarray*}
\sho_X^\lambda
=\rhom[\shd_{\overline{X}}](\sho_{\overline{X}},\shc^{\infty,\lambda}_X),
\,\,\lambda=t,{\rm w},\omega.
\end{eqnarray*}
Then, $\sho_X^t,\sho_X^{\rm w}$ and $\sho_X^\omega$ belong to 
$D^b(\md[\shd_X,{\II[\C_X]}])$.
Note that if $F$ is $\R$-constructible, 
\begin{eqnarray*}
\rhom(F,\sho_X^t)&\simeq& \sht\hom(F,\sho_X)\\
\rhom(F,\sho_X^{\rm w})&\simeq&D'F\wtens\sho_X\\
\rhom(F,\sho_X^\omega)&\simeq&D'F\tens \sho_X,
\end{eqnarray*}
where $\sht\hom(\cdot,\sho_X)$ and $\cdot\wtens\sho_X$ are the 
functors of
tempered  and formal cohomology of \cite{K1} and \cite{K-S2}, respectively.

In particular, let $M$ be  a real analytic manifold and assume that $X$
is a complexification of $M$.
We find:
\begin{eqnarray*}
\rhom(D'\C_M,\sho_X)&\simeq& \shb_M 
\mbox{ (Sato's hyperfunctions)},\\
\rhom(D'\C_M,\sho_X^t)&\simeq& \shd b_M
\mbox{ (Schwartz's distributions)},\\
\rhom(D'\C_M,\sho_X^{\rm w})&\simeq&\shc^\infty_M
\mbox{ ($C^\infty$-functions)},\\
\rhom(D'\C_M,\sho_X^\omega)&\simeq&\sha_M
\mbox{ (real analytic functions)}.
\end{eqnarray*}
Replacing $\hom$ by $\ihom$ we find new ind-sheaves.
For example, we can define the ind-sheaf $\shd b^t_M$ 
of tempered distributions by setting
$$\shd b^t_M={\rm R}\ihom(D'\C_M,\sho_X^t).$$
Note that $\alpha_X(\shd b^t_M)\simeq \shd b_M$.

\begin{remark}
The object $\sho_X^t$ is not concentrated in degree $0$ if $\dim X>1$.
In fact if $\shc$ is an abelian category, a complex 
$A'\to[f] A\to[g] A''$ in $\indc$ is exact 
if and only if the dotted arrows in the diagram below with $X\in \shc$
may be completed with $Y\in \shc$ in such a way that the morphism
$\alpha:Y\to X$ is an epimorphism:
$$\xymatrix{
A'\ar[r]^f &A\ar[r]^g&A''\\
Y\ar@{.>}[u]^t\ar@{.>>}[r]_\alpha&X\ar[u]^s\ar[ru]|0&
}$$
Let us apply this result to the complex
$$\shc^{\infty,t,(p-1)}_X \to[\overline{\partial}]
\shc^{\infty,t,(p)}_X \to[\overline{\partial}]\shc^{\infty,t,(p+1)}_X
$$
and choose $F=\C_U$, $U$ open in $X$. Assume $\sho_X^t$ is
concentrated in degree $0$. Then, for any $s\in  \shc^{\infty,t,(p)}(U)$
solution of $\overline{\partial}s=0$, there exists an epimorphism 
$\alpha:G\epito F$
and $t\in\Hom(G,\shc^{\infty,t,(p-1)})$ such that 
$s\circ \alpha=\overline{\partial}t$.
We may assume $G$ is a finite direct sum of sheaves $\C_{U_j}$, with
$U=\cup_jU_j$. We thus find that 
$s\vert_{U_j}=\overline{\partial}t_j$, which is in general not possible.

Note that the same argument hold with the sheaf $\sho_X$, which shows
that the functor $\rho$ in (\ref{functrho}) is not exact.
\end{remark}

\section{Microlocalization}

In this section, $X$ denotes a real analytic manifold.
We denote by $p_1$ and $p_2$ the first and second projection defined
on $X\times X$ and by $\Delta$ the diagonal. 
We denote by $\tau:TX\to X$
and $\pi:T^*X\to X$
its tangent and cotangent bundles, respectively. We denote by $T_YX$
and $T^*_YX$ the normal and conormal bundle to a closed submanifold $Y$
of $X$. In particular, $T_XX$ and $T^*_XX$ are the 
zero-sections of these bundles. We shall identify $T_\Delta(X\times X)$ 
with $TX$ by the first projection on $TX\times TX$.
If $S$ is a subset of $X$, one denotes by 
$C_Y(S)$ the Whitney normal cone of $S$ along $Y$, a closed cone in $T_YX$.

We denote by $\omega_X$
the dualizing complex on $X$. (Hence, $\omega_X\simeq or_X[\dim X]$.)
The micro-support $SS(F)$ of a sheaf $F$ and the functor $\mu hom$ are
defined in \cite{K-S1}.

For a section $s:X\to T^*X$
of $\pi$, one can construct an object $L_s\in D^b(\II[k_{X\times X}])$
with the properties that $\supp L_s=\Delta$ and such that 
\begin{eqnarray*}
&&L_s\otimes \widetilde{\C}_{s^{-1}(T^*_XX)}\simeq
\widetilde{\C}_{s^{-1}(T^*_XX)}\\
&&L_s\simeq \inddlim[U,V]k_{U\cap\overline{V}}\tens\opb{p_2}\omega_X
\quad\mbox{over $X\setminus {s^{-1}(T^*_XX)}$},
\end{eqnarray*}
where $V$ ranges through the family of open neighborhoods of $\Delta$
and $U$ through the family of open subsets of $X\times X$ such that:

$$C_\Delta(U)\subset \{v\in T_xX;\langle v,s(x)\rangle<0\}\cup T_XX.$$
 One defines 
$$L_s\circ G=R\eeim{p_1}(L_s\tens\opb{p_2}G),$$
and one proves that for $F,G\in\II[k_X]$ 
$$\opb{s}\mu hom(F,G)\simeq \rhom(F,L_s\circ G).$$
Now assume that $s:T^*X\to T^*(T^*X)$ is a section. Then 
$\opb{s}\mu hom(\opb{\pi}F,\opb{\pi}G)\simeq \mu hom(F,G)$.
Taking as $s$ the canonical $1$-form $\alpha_X$ on $T^*X$, one can then
construct an object $K_X\in D^b(\II[k_{T^*X\times T^*X}])$ and define 
the microlocalization functor
\begin{eqnarray*}
\mu_X:D^b(\II[k_X])\to D^b(\II[k_{T^*X}]),\\
\mu_X(F)=K_X\circ \opb{\pi}F.
\end{eqnarray*}

\begin{theorem}
Let $F,G\in D^b(k_X)$. Then
\begin{eqnarray*}
SS(F)&=&\supp\mu_X(F),\\
\mu hom(F,G)&\simeq& \rhom(\mu_X(F),\mu_X(G)),\\
              &\simeq& \rhom(\opb{\pi}F,\mu_X(G)).
\end{eqnarray*} 
\end{theorem}

\section{Applications}
Let $X$ denote  a complex manifold of complex dimension $n$. 
On $T^*X$ there are some classical sheaves
associated with a sheaf $F$ on $X$: the sheaf
$\mu hom(F,\sho_X)$, or ($F$ being $\R$-constructible) the sheaves
$t\mu hom(F,\sho_X)$ of \cite{An} and $w\mu hom(F,\sho_X)$ of \cite{Co}.
We can obtain all these sheaves in a unified way:
\begin{eqnarray*}
\mu hom(F,\sho_X)
&\simeq&\rhom(\opb{\pi}F,\mu_X(\sho_X)),\\
t\mu hom(F,\sho_X)
&\simeq&\rhom(\opb{\pi}F,\mu_X(\sho^t_X)),\\
w\mu hom(F,\sho_X)
&\simeq&\rhom(\opb{\pi}F,\mu_X(\sho^{\rm w}_X)).
\end{eqnarray*}
In particular, the sheaf $\she_X^\R$ of microlocal operators 
of \cite{S-K-K} is isomorphic to the sheaf
$\rhom(\opb{\pi}\C_\Delta,\mu_{X\times X}(\sho_{X\times X}^{(0,n)}))[n]$.
\begin{theorem}
The complex $\mu_X(\sho_X)[n]$ is concentrated in degree $0$ on
$\dot{T}^*X:=T^*X\setminus T^*_XX$. Moreover, 
$\mu_X(\sho_X)[n]\in \md[{\she_X^\R,\II[\C_{\dot{T}^*X}]}]$.
\end{theorem}
\begin{corollary}
Let $F\in D^b(\C_X)$. The object $\mu hom(F,\sho_X)$ is well-defined
in the derived category $D^b(\she_X^\R)$.
\end{corollary}

\ifx\undefined\bysame
\newcommand{\bysame}{\leavevmode\hbox to3em{\hrulefill}\thinspace}
\fi

\bigskip
\bigskip
\bigskip

\noindent {\small M-K Research Institute for Mathematical Sciences\\
Kyoto University, Kyoto 606 Japan}

\bigskip

\noindent {\small P-S Institut de Math\'ematiques, Analyse Alg\'ebrique\\
Universit\'e P \& M Curie, Case 82\\4, place Jussieu F-75252, Paris
Cedex 05, France

\noindent {\tt schapira@math.jussieu.fr \ \ 
 http://www.math.jussieu.fr/\~{}schapira/} }

\end{document}